\documentclass[12pt]{article}
\usepackage{amsmath}
\usepackage{amsfonts}
\usepackage{amssymb}
\usepackage{theorem}

\title{Nice enumerations and extreme amenability}
\author{A. Ivanov
\thanks{The author is supported by Polish National Science Centre grant DEC2011/01/B/ST1/01406} } 
\date{ } 

\newtheorem{thm}{Theorem}[section] 
\newtheorem{lem}[thm]{Lemma}
\newtheorem{definicja}[thm]{Definition}

\newtheorem{prop}[thm]{Proposition} 
\newtheorem{remark}[thm]{Remark}

\begin{document}
\maketitle
\topskip 20pt

\begin{quote}
{\bf Abstract.} 
We study properties related to nice enumerability of countably 
categorical structures  and properties related to extreme amenability 
of automorphism groups of  these structures.   

\bigskip


\end{quote}

\bigskip


\section{Introduction} 

A group $G$ is called {\bf amenable} if every $G$-flow 
(i.e. a compact Hausdorff space along with a continuous G-action) 
supports an invariant Borel probability measure.  
If  every $G$-flow has a fixed point then we say that $G$ 
is {\bf extremely amenable}. 
Let $M$ be a relational countably categorical structure which is 
a Fra\"{i}ss\'{e} limit of a Fra\"{i}ss\'{e} class $\mathcal{K}$. 
In particular $\mathcal{K}$ coincides with $Age(M)$, 
the class of all finite substructures of $M$. 
By Theorem 4.8 of the paper of Kechris, Pestov and Todorcevic \cite{KPT} 
the group $Aut(M)$ is {\em extremely amenable if and only if 
the class $\mathcal{K}$ has the Ramsey property and consists of rigid elements.} 
Here the class $\mathcal{K}$ is said to have 
the {\bf Ramsey property} if  if for any $k$ 
and a pair $A<B$ from $\mathcal{K}$ 
there exists $C\in \mathcal{K}$ so that 
each $k$-coloring 
$$
\xi :{C\choose A}\rightarrow k
$$ 
is monochromatic on some ${B'\choose A'}$
from $C$ which is a copy of 
${B\choose A}$, i.e. 
$$ 
C\rightarrow (B)^{A}_k . 
$$ 
\parskip0pt 

We remind the reader that a $G$-flow  
$X$ is called {\bf minimal}, if every its $G$-orbit is dense. 
The flow $X$ is {\bf universal}, if for every $G$-flow $Y$ 
there is a continuous $G$-map $f: X \rightarrow Y$. 
According to topological dynamics a universal minimal 
flow always exists and is unique up to $G$-flow isomorphism 
(and is usually denoted by $M(G)$). 
The following question was formulated by several people. 
In particular it appears in the paper of 
Angel, Kechris and Lyons \cite{AKL}.  
\begin{quote} 
Let $G=Aut(M)$, where 
$M$ is a countably categorical structure. 
Is the universal minimal $G$-flow  metrizable?  
\end{quote}  

Recently A.Zucker has found a characterisation 
of automorphism groups of relational structures 
which have metrizable universal minimal flow. 
It substantially develops the previous work 
of Kechris, Pestov, Todorcevic and Nguyen van Th\'{e} 
from \cite{KPT} and  \cite{The}. 
 
\bigskip 

{\bf Theorem A} (Theorem 1.2 of \cite{ZA}).  
{\em Let $M$ be a relational structure which is 
a Fra\"{i}ss\'{e} limit of a Fra\"{i}ss\'{e} class $\mathcal{K}$. 
Then the following are equivalent.  

1) $G=Aut(M)$ has mertizable universal minimal flow, 

2) 
there is a sequence of new relational symbols 
$\bar{S}$ and a precompact  $\bar{S}$-expansion of $M$, 
say $M^*$,  so that 
\begin{quote} 
(i) $M^*$ is a Fra\"{i}ss\'{e} structure, 

(ii) $Aut(M^* )$ is extemely amenable  and  

(iii) the closure of the $G$-orbit 
of $M^*$ in the space of $\bar{S}$-expansions of $M$ 
is a universal minimal $G$-flow. 
\end{quote} 
Moreover if $M(G)$ is metrizable, then $G$ 
has the {\bf generic point property}, i.e.  
$M(G)$ has a $G_{\delta}$-orbit. }    
\bigskip 

In this formulation precompactness means that every 
member of $\mathcal{K}$ has finitely many expansions 
in $Age(M^* )$  

By this theorem it is crucial to know whether there 
is a countably categorical structure $M$ which 
does not have expansions as in Theorem A. 
It is worth noting that some versions of this question were 
formulated for example  in \cite{BPT}, see Problems 27, 28.   
Related results can be also found in \cite{KS},  \cite{AKL} 
and \cite{Z}.  
\parskip0pt 

In our paper having in mind these respects, 
we consider automorphism groups of countably 
categorical structures which satisfy some properties 
related to nice enumerability, \cite{az}.

\begin{definicja} 
Let $M$ be a countable structure. 
A linear ordering $\prec$ of $M$ of type $\omega$ 
is called an {\bf AZ-enumeration} of $M$ 
if for any $n\ge 1$ it satisfies the following property:
\begin{quote}
whenever $\bar{b}_i$, $i<\omega$, is 
a sequence of $n$-tuples from $M$,  
there exist some $i<j<\omega$ and a $\prec$-preserving 
$Th(M)$-elementary map $f:M \rightarrow M$
such that $f(\bar{b}_i )=\bar{b}_j$.
\end{quote}  
\end{definicja}

It is known that any $\omega$-dimensional 
classical geometry has an AZ-enumeration, 
which is obtained by some canonical procedure, 
see \cite{CH}. 
Moreover the class of finite dimensional geometries 
of a fixed type considered with canonical orderings 
is a Fra\"{i}ss\'{e} class with the property 
that canonical orderings are unique up to isomorphism. 
In the case of vector spaces over a finite field  
this implies the Ramsey property of the class 
with canonical orderings \cite{KPT}. 
This example was the starting point of the paper. 
Is AZ-enumerability connected with the Ramsey property? 
We will see below that at least these properties 
have similar consequences. 
In fact they imply some kind of amalgamation 
of ordered expansions. 
This will be proved in Section 2. 

In Section 3 we study a construction which 
produces orderings with properties similar to AZ-enumerability.  
Since the construction uses some amalgamation, 
the final orderings are not of type $\omega$. 
We introduce more general {\em AZ-quasi-well-orderings} 
and show in Section 3.3 how they can be applied 
for permutation modules. 
Our motivation here is as follows. 
When a countably categorical structure 
has a nice enumeration (in particular when 
it has an AZ-enumeration), it can be studied 
by some standard methods, see \cite{az}, \cite{CaEv}, \cite{CH}. 
Since it is not known if any countably 
categorical structure  has a nice enumeration, 
we try to find a weaker condition, 
which is some kind of 
{\em potential enumerability}. 
Properties related to extreme amenability 
are helpful in this approach. 

In Section 1 we describe some methods 
which will be used in the paper. 
Basically they are taken from \cite{iv}. 
Then we slightly modify the approach to extreme 
amenability from \cite{KPT} and \cite{The}  
so that it works for expansions 
of structures where elimination of quantifiers 
is not necessarily satisfied, 
for example obtained by Hrushovski's 
amalgamation method.  
This brings additional flexibility. 
Here we also use \cite{iv} 
(and of course \cite{kechros}).

\section{Generic expansions of $\omega$-categorical 
structures and extreme amenability}

We fix a countable structure $M$ in a language $L$.
We assume that $M$ is $\omega$-categorical
(most of the terms below make sense under the assumption that
$M$ is atomic).
Let $T$ be an extension of $Th(M)$ in the language with additional
relational and functional symbols 
${\bf \bar{r}} = ({\bf r}_{1},...,{\bf r}_{t})$.
We assume that $T$ is axiomatizable
by sentences of the following form:
$$
(\forall \bar{x})(\bigvee_{i} (\phi_{i}(\bar{x})\wedge\psi_{i}(\bar{x}))),
$$
where $\phi_{i}$ is a quantifier-free formula in the language
$L\cup {\bf \bar{r}}$, and $\psi_{i}$ is a first-order formula of
the language $L$. 
Consider the set ${\bf X}$ of all possible expansions
of $M$ to models of $T$. \parskip0pt

Following \cite{iv} we define for a tuple $\bar{a}\subset M$
a {\em diagram} $\phi(\bar{a})$ of ${\bf \bar{r}}$ on $\bar{a}$.
To every functional symbol from ${\bf \bar{r}}$ we associate
a partial function from $\bar{a}$ to $\bar{a}$.
Choose a formula from every pair
$\{ {\bf r}_{i}(\bar{a}'),\neg {\bf r}_{i}(\bar{a}')\}$, where 
${\bf r}_{i}$ is a relational symbol from ${\bf \bar{r}}$ and 
$\bar{a}'$ is a tuple from $\bar{a}$ of the corresponding length.
Then $\phi(\bar{a})$ consists of the conjunction of the chosen
formulas and the definition of the chosen functions
(so, in the functional case we look at $\phi(\bar{a})$ as
a tuple of partial maps). \parskip0pt

Consider the class ${\bf B}_{T}$ of all theories
$D(\bar{a}), \bar{a}\subset M,$ such that each of them
consists of $Th(M,\bar{a})$ and a diagram of ${\bf \bar{r}}$ 
on $\bar{a}$ satisfied in some $(M,{\bf \bar{r}})\models T$.
We order ${\bf B}_{T}$ by extension: $D(\bar{a})\le D'(\bar{b})$
if $\bar{a}\subset \bar{b}$ and $D'(\bar{b})$ implies
$D(\bar{a})$ under $T$ (in particular, the partial functions
defined in $D'$ extend the corresponding partial functions
defined in $D$).
Since $M$ is an atomic model, each element of ${\bf B}_{T}$ is
determined by a formula of the form
$\phi(\bar{a})\wedge\psi(\bar{a})$, where $\psi$ is a complete
formula for $M$ and $\phi$ is a diagram of ${\bf \bar{r}}$ on $\bar{a}$.
The corresponding formula $\phi(\bar{x})\wedge\psi(\bar{x})$
will be called {\em basic}. \parskip0pt

On the set ${\bf X} =\{ (M,{\bf \bar{r}'}):(M,{\bf \bar{r}'})\models T \}$
of all ${\bf \bar{r}}$-expansions of the structure $M$ we consider
the topology generated by basic open sets of the form 
$$ 
[D(\bar{a})]=\{ (M,{\bf \bar{r}'}):(M,{\bf \bar{r}'})\models D(\bar{a})\} 
\mbox{ , }  \bar{a}\subset M.
$$ 
It is easily seen that any $[D(\bar{a})]$ is clopen.
The topology is metrizable: fix an enumeration
$\bar{a}_{0},\bar{a}_{1},...$ of $M^{< \omega}$ and define
\begin{quote}
$d((M,{\bf \bar{r}'}),(M,{\bf \bar{r}''}))=\sum \{ 2^{-n}:$
there is a symbol ${\bf r}\in {\bf \bar{r}}$ such that its 
interpretations on $\bar{a}_{n}$ in the structures $(M,{\bf \bar{r}'})$ 
and $(M,{\bf \bar{r}''})$ are not the same (if ${\bf r}$ is a functional 
symbol then ${\bf r'}(\bar{b}) \not= {\bf r''}(\bar{b})$ for some
$\bar{b}\subseteq \bar{a}_{n}$) $\}$.
\end{quote}

\noindent
It is easily seen that the metric $d$ defines the topology
determined by the sets of the form $[D(\bar{a})]$. \parskip0pt

By the assumptions on $T$ ($T$ is axiomatizable by sentences
which are universal with respect to symbols from ${\bf \bar{r}}$)
the space ${\bf X}$ forms a closed subset of the complete metric 
space of all ${\bf \bar{r}}$-expansions of $M$. 
Thus ${\bf X}$ is complete and the Baire Category
Theorem holds for ${\bf X}$. 
We say that $(M,{\bf \bar{r}})\in {\bf X}$ is {\it generic} if
the class of its images under $Aut(M)$ is comeagre in ${\bf X}$ 
\cite{iv}.\parskip0pt

\begin{remark} \label{precomp} 
{\em All our arguments also work for the case when 
${\bf \bar{r}} = ({\bf r}_{1},...,{\bf r}_{t},...)$
is an infinite sequence, but for every tuple $\bar{b}$ 
from $M$ the family ${\bf B}_T$ has finitely many 
diagrams defined on $\bar{b}$. 
In this case we say that ${\bf X}$ consists of 
{\em precompact} expansions. 
} 
\end{remark} 

Notice that the space $Aut(M)$ under the conjugacy action, and 
generic automorpisms (introduced in \cite{truss}) provide 
a particular example of this construction.
Indeed, identify each $\alpha \in Aut(M)$ with the expansion
$(M,\alpha,\alpha^{-1})$.
The class of structures of this form is
axiomatized in the language of $M$ with the functional symbols
$\{ \alpha, \beta \}$ by $Th(M)$, the sentence
$\alpha \beta(x) = \beta \alpha(x) = x$ and universal sentences
asserting that $\alpha$ preserves the relations of $M$.
Also, any partial isomorphism $\bar{a}\rightarrow\bar{a}'$
can be viewed as the diagram corresponding to the maps
$\bar{a}\rightarrow\bar{a}'$ and $\bar{a}'\rightarrow\bar{a}$.
It is clear that a generic automorphism $\alpha$ (see \cite{truss}) 
defines generic expansion $(M,\alpha,\alpha^{-1})$.
\parskip0pt

Similar considerations can be applied in the following general 
situation.
Let $M$ be an $\omega$-categorical structure in a language $L$.
Let ${\bf \bar{r}}$ be a tuple of relations on $M$ and $T$ be $Th(M)$
extended by all the sentences from $Th(M,{\bf \bar{r}})$ of the form
$\forall \bar{x} \neg D(\bar{x})$, where $D(\bar{x})$ is basic
for $(M,{\bf \bar{r}})$.
It is clear that $T$ satisfies the conditions from
the beginning of the section.
Note that ${\bf B}_{T}$ consists of all diagrams $D(\bar{b})$
such that the corresponding formula $D(\bar{x})$ is realizable
in $(M,{\bf \bar{r}})$.
The expansion $(M,{\bf \bar{r}})$ is {\it ubiquitous in category} if
$(M,{\bf \bar{r}})$ is generic with respect to ${\bf B}_{T}$ \cite{iv}.
It is clear that generic expansions (with respect to some theory $T'$) 
are always ubiquitous in category. \parskip0pt

Theorem 1.5 of \cite{iv} states that
\begin{quote}
a structure $(M,{\bf \bar{r}})$ is ubiquitous in category 
if and only if every complete type over $\emptyset$ 
realizable in $(M,{\bf \bar{r}})$ is determined in 
$(M,{\bf \bar{r}})$ by a formula of the form 
$\exists \bar{y} D(\bar{x}\bar{y})$ where 
$D(\bar{x}\bar{y})$ is basic.
\end{quote}

We now give several important definitions from \cite{iv}. 
We say that ${\bf B}_{T}$ has the {\em joint embedding property} 
if for any $D_{1}(\bar{a}), D_{2}(\bar{b}) \in {\bf B}_{T}$ 
there exist $D(\bar{c}) \in {\bf B}_{T}$ and $M$-elementary 
maps $\delta: \bar{a} \rightarrow \bar{c}$ and
$\sigma: \bar{b} \rightarrow \bar{c}$ such that $D(\bar{c})$
extends $D_{1}(\delta(\bar{a})) \cup D_{2}(\sigma(\bar{b}))$. 
The class ${\bf B}_{T}$ has the {\em weak amalgamation property} 
(see \cite{kechros}, in the original paper \cite{iv} it is called 
the {\em almost amalgamation property}) if for every 
$D(\bar{a}) \in {\bf B}_{T}$ there is an extension
$D'(\bar{a}\bar{b}) \in {\bf B}_{T}$ such that for any
$D_{1}(\bar{a}\bar{c}_{1}), D_{2}(\bar{a}\bar{c}_{2}) \in {\bf B}_{T}$, 
where $D'(\bar{a}\bar{b}) \le D_{i}(\bar{a}\bar{c}_{i}), i = 1,2$, 
there exists a common extension $D''(\bar{a}\bar{c}) \in {\bf B}_{T}$ 
under some $(M,\bar{a})$-elementary maps 
$\bar{c}_{i} \rightarrow \bar{c}, i = 1,2$ (which may move $\bar{b}$).

\bigskip

{\bf Theorem B.} (\cite{iv}, Theorem 1.2 and Corollary 1.4)
{\em 
(a)  The set ${\bf X}$ has a generic structure if and
only if ${\bf B}_{T}$ has the joint embedding property 
and the weak amalgamation property.

(b) If there are no continuum many pairwise non-isomorphic elements 
of ${\bf X}$, then ${\bf B}_T$ has the weak amalgamation property. }

\bigskip 

It is worth noting that in (a) for any member of ${\bf B}_T$ 
the corresponding basic formula is realised in a generic 
structure.

An element $D(\bar{b}) \in {\bf B}_{T}$ is called an 
{\em amalgamation base} if any two of its extensions 
have a common extension in ${\bf B}_{T}$ under some 
automorphism of $M$ fixing $\bar{b}$. 
We say that ${\bf B}_{T}$ satisfies {\em Truss' condition}  
if any element of ${\bf B}_{T}$ extends to 
an amalgamation base. 
If it holds then the set of amalgamation bases is 
a cofinal subset of ${\bf B}_{T}$ which has 
the amalgamation property. 
It is clear that Truss' condition implies 
the weak amalgamation property. 
In particular it together with the joint 
embedding property implies the existence 
of a generic expansion $(M,\bar{{\bf r}})$. 
The proof of Theorem 1.5 of \cite{iv} 
shows that when ${\bf B}_{T}$ satisfies Truss'  
condition and $D(\bar{b})\in {\bf B}_{T}$ is an amalgamation base, 
then the type of $\bar{b}$ in this $(M,\bar{{\bf r}})$
is determined by the basic formula $D(\bar{x})$. 
We will use this fact later.  

The group $Aut(M)$ has a natural action on 
${\bf B}_T$. 
Moreover if $\bar{a}$ and $\bar{b}$ have the same 
type with respect to $Th(M)$, then 
after replacement $\bar{a}$ by $\bar{b}$ 
in $D(\bar{a})$ we obtain an image of 
$D(\bar{a})$ under an automorphism 
of $M$. 
This is a consequence of an $\omega$-homogeneity 
of $M$. 
It is also clear that $Aut(M)$ acts continuously 
on ${\bf X}$ with respect to the topology 
defined above.  

Let $D(\bar{a})\in {\bf B}_T$ have 
an extension $D'(\bar{b})\in {\bf B}_T$. 
Let ${D'(\bar{b})\choose D(\bar{a})}$ 
be the set of all images of $D(\bar{a})$
in $D'(\bar{b})$ under elementary maps 
of $M$. 

\begin{definicja} \label{Ramsey} 
An $Aut(M)$-invariant  subfamily $\mathcal{C} \subseteq {\bf B}_T$ 
satisfies the Ramsey property if for any $k$ 
and a pair $D(\bar{a})< D'(\bar{b})$ from $\mathcal{C}$ 
there exists $D''(\bar{c})\in \mathcal{C}$ so that 
each $k$-coloring 
$$
\xi :{D''(\bar{c})\choose D(\bar{a})}\rightarrow k
$$ 
is monochromatic on some ${D'(\bar{b}')\choose D(\bar{a}')}$
from $D''(\bar{c})$ which is a copy of 
${D'(\bar{b})\choose D(\bar{a})}$, i.e. 
$$ 
D''(\bar{c})\rightarrow (D'(\bar{b}))^{D(\bar{a})}_k . 
$$ 
\end{definicja} 

The following lemma follows from the 
argument used in the proof of 
Theorem 1.2 ($1 \rightarrow 2$) of 
\cite{HubNes} 
(see also Proposition 8.13 of \cite{ZA}). 

\begin{lem} \label{HuNe} 
Assume that an $Aut(M)$-invariant  subfamily 
$\mathcal{C} \subseteq {\bf B}_T$ 
satisfies the joint embedding property and 
the Ramsey property.  
Then $\mathcal{C}$ 
satisfies the amalgamation property.  
\end{lem} 

\begin{remark} \label{R->A}
{\em By the proof of Theorem 1.2 ($1 \rightarrow 2$) 
of \cite{HubNes} the weak amalgamation property 
is a consequence of the joint embedding property 
and following version of 
the Ramsey property: 
\begin{quote} 
For any $D(\bar{a})\in {\bf B}_T$ there is 
an extension $D'(\bar{a}\bar{b}) \in {\bf B}_{T}$ 
such that for any $D_{1}(\bar{a}\bar{c}_{1}) \in {\bf B}_{T}$, 
where $D'(\bar{a}\bar{b}) \le D_{1}(\bar{a}\bar{c}_{1})$, 
there exists an extension 
$D_1 (\bar{a}\bar{c}_1 )< D_2(\bar{a}\bar{c}_2) \in {\bf B}_{T}$ 
such that 
$$ 
D_2 (\bar{a}\bar{c}_2 )\rightarrow ( D_1 (\bar{a}\bar{c}_1 ))^{D(\bar{a}) }_2 .
$$ 
under some $(M,\bar{a})$-elementary maps. 
\end{quote} 
}
\end{remark} 

The following theorem is a slightly generalized  
version of Theorem 4.5 from \cite{KPT}. 

\begin{thm} \label{extram} 
Let $M$ be an $\omega$-categorical structure 
and ${\bf B}_T$ satisfy Truss' condition. 
Let $\mathcal{C}\subset {\bf B}_T$ be an invariant  
cofinal subset of  amalgamation bases with the joint 
embedding property and the amalgamation property.  

Then the automorphism group $Aut(M,\bar{{\bf r}})$ 
of a generic expansion corresponding to 
$\mathcal{C}$ is extremely amenable if and only 
if the class  $\mathcal{C}$ has the Ramsey property 
and consists of rigid elements, i.e. no 
$D(\bar{a})\in \mathcal{C}$ can be taken onto itself 
by a non-trivial elementary map $\bar{a}\rightarrow \bar{a}$ 
with respect to $Th(M,\bar{{\bf r}})$. 
\end{thm} 

{\em Proof.} 
We adapt the proof of Theorem 4.5 from \cite{KPT} as follows. 
Firstly we remind the reader that for a closed subgroup $G<Aut(M)$ 
a $G$-{\em type} of a tuple $\bar{a}$ is just the orbit $G\bar{a}$.   
It is a consequence of the proof of Theorem 1.5 from \cite{iv} 
that for $G=Aut(M,\bar{{\bf r}})$ if  $D(\bar{a}) \in \mathcal{C}$, 
then  the $G$-type of $\bar{a}$ is determined by $D(\bar{x})$, 
i.e. coincides with the set of all realisations of this formula in 
$(M,\bar{{\bf r}})$. 
Then the ordering of $G$-types $G\bar{a} \le G\bar{b}$ 
introduced in \cite{KPT} just corresponds to the relation 
$D(\bar{y}) \vdash D(\bar{x})$ for some embedding 
of $\bar{x}$ (corresponding to $\bar{a}$) into 
$\bar{y}$ (corresponding to $\bar{b}$). 
The Ramsey property introduced in \cite{KPT} 
in the case of $Aut(M,\bar{{\bf r}})$-types of 
tuples $\bar{a}$ with $D(\bar{a})\in \mathcal{C}$ 
coincides with Definition \ref{Ramsey} for $\mathcal{C}$. 
The rest follows from Proposition 4.3, Theorem 4.5 and Remark 4.6 
of \cite{KPT}. 
$\Box$ 

\bigskip

It is worth noting that the majority of 
basic statements of \cite{KPT} 
can be adapted to the situation of 
Theorem \ref{extram}. 
In particular Theorem 7.5 of \cite{KPT} 
and Theorem 5 of \cite{The} 
(a generalisation of the former one)   
can be stated as follows. 

\begin{thm} \label{companion} 
Let $M$ be an $\omega$-categorical structure, 
$G=Aut(M)$ and ${\bf B}_T$ satisfy Truss' condition. 
Let $\mathcal{C}\subset {\bf B}_T$ be an invariant  
cofinal subset of  rigid amalgamation bases with the joint 
embedding property and the amalgamation property.  
Assume that $(M,\bar{{\bf r}})$ is a generic expansion 
corresponding to $\mathcal{C}$. 

Then the space ${\bf X}$ is the universal minimal flow 
of  $G$  if and only if the class  $\mathcal{C}$ has 
the Ramsey property and the following {\bf expansion property} 
(relative to $M$): 
\begin{quote}
any tuple $\bar{a}$ from $M$ extends to a tuple $\bar{b}\in M$ 
so that for any $D(\bar{a})$ and $D'(\bar{b}) \in \mathcal{C}$ 
there is an $M$-elementary map $\bar{a}\rightarrow \bar{b}$ 
which embeds $D(\bar{a})$ into $D'(\bar{b})$. 
\end{quote} 
\end{thm} 

A proof of this theorem can be obtained by a straightforward 
adaptation of Sections 4 and 5 of \cite{The}. 
Moreover repeating Proposition 14.3 from \cite{AKL} 
we can show that in the situation of Theorem \ref{companion} 
the group $G$ has the generic point property, i.e.  
every minimal $G$-flow has a comeagre orbit. 

\begin{remark} 
{\em 
Theorem \ref{companion} can be also considered 
as a generalisation of Theorem 5.7 of \cite{ZA}. 
To see this one needs some routine work to connect 
our terminology  with that from \cite{ZA}. 
In particular a class $\mathcal{K}$ of finite structures 
is {\em Fra\'{i}ss\'{e} - HP} in terms of \cite{ZA} 
iff there is an age $\mathcal{K}'$ so that 
$\mathcal{K}$ is a cofinal subclass of amalgamation 
bases in $\mathcal{K}'$. 
In particular $\mathcal{K}'$ satisfies Truss' condition.  
}
\end{remark}

\section{Generic expansions} 

We now generalize the definition of AZ-enumerations. 
The previous one appears exactly as the partial case 
described in this definition.  

\begin{definicja} 
Let $M$ be a countable structure. 
A linear ordering $\prec$ of a countable subset 
$\mathcal{L}\subseteq M$ is called an {\bf AZ-ordering} in $M$ 
if for any $n\ge 1$ it satisfies the following property:
\begin{quote}
whenever $\bar{b}_i$, $i<\omega$, is a sequence of
$n$-tuples from $\mathcal{L}$, which is not decreasing coordinatewise, 
there exist some $i<j<\omega$ and a $\prec$-preserving 
$Th(M)$-elementary map $f:\mathcal{L} \rightarrow \mathcal{L}$
such that $f(\bar{b}_i )=\bar{b}_j$.
\end{quote} 
When additionally $M=\mathcal{L}$ and $\prec$ has order-type 
$\omega$ we say that $\prec$ is an {\bf AZ-enumeration} of $M$.  
\end{definicja}

\subsection{AZ-enumerations and generic expansions}

The following theorem shows that AZ-enumerations 
provide generic expansions by orderings. 
In particular we have a version of the amalgamation 
property for expansions (WAP). 

\begin{thm} \label{GenStr} 
Let $M$ be a Fraisse limit of a relational 
amalgamation class $\mathcal{K}$. 
Let $\prec$ be a linear ordering of some $I\subseteq M$. 
Let  $T$ be $Th(M)$ extended by all the 
sentences from $Th(M, \prec )$ of the form 
$\forall \bar{x} \neg D(\bar{x})$, where 
$D(\bar{x})$ is basic for $(M,\prec )$.

If the linear ordering $(I, \prec )$ has the property that 
\begin{quote}
whenever $\bar{a}_i$, $i<\omega$, is a sequence of
$n$-tuples from $M$, there exist some $i<j<\omega$ and
an  elementary map $f:M\rightarrow M$ which is an isomorphic 
injection of $(M,I,\prec )$ into itself 
such that $f(\bar{a}_i )=\bar{a}_j$, 
\end{quote}
then  the space 
of all $T$-expansions of $M$ contains a generic 
structure which is $\omega$-categorical. 

In particular this conclusion holds when $\prec$ 
is an AZ-enumeration of $M$. 
\end{thm} 

{\em Proof.} 
Assume that $\prec$ is an ordering of $I\subseteq M$ as in the formulation. 
The definition of $T$ implies that ${\bf B}_T$ has JEP. 
Let us show that ${\bf B}_T$ satisfies WAP, 
the weak amalgamation property. 
If WAP does not hold for some $D(\bar{a})\in {\bf B}_T$, 
then we build a tree in $({\bf B}_T ,<)$ 
with the root $D(\bar{a})$. 
At every step we split the already constructed 
extensions of $D(\bar{a})$ as follows.  
If $D'(\bar{a}\bar{b})$ is an extension 
corresponding to the vertex 
$\varepsilon_1 \varepsilon_2 ...\varepsilon_n$
of the tree (with $\varepsilon_i \in \{ 0,1\}$), 
find $D_1(\bar{a}\bar{c}_1)$ and $D_2(\bar{a}\bar{c}_2)$ 
extending $D'(\bar{a}\bar{b})$ which cannot be 
amalgamated over $D(\bar{a})$. 
These extensions correspond to vertices 
$\varepsilon_1 \varepsilon_2 ...\varepsilon_n 0$ 
and 
$\varepsilon_1 \varepsilon_2 ...\varepsilon_n 1$. 

Now choose tuples $\bar{a}_i$ from $(M,I, \prec )$ 
which correspond to $\bar{a}$ in all 
extensions $D'(\bar{a}\bar{b})$ with numbers 
$0$, $10$, ..., $111...10$, ... in the tree. 
In other words $\bar{a}_i$ extends to 
a tuple $\bar{a}_i \bar{b}'$ realizing 
$D'(\bar{x}\bar{y})$ corresponding 
to the number $1...10$ with $i$ units. 
It is clear that no $\bar{a}_i$ can be 
taken to $\bar{a}_j$ with $i<j$ by 
an  elementary map, which is 
an isomorphic embedding of $(M, I, \prec )$
into itself. 
This is a contradiction with the assumptions. 

By Theorem B find a generic structure $(M,I^* ,\prec^* )$. 
If this structure is not 
$\omega$-categorical, then 
for some natural $k>0$ there are infinitely 
many $k$-types over $\emptyset$. 
By Theorem 1.5 of \cite{iv} each type over 
$\emptyset$ realizable in $(M,I^* ,\prec^* )$ 
is determined by a formula of the form 
$\exists \bar{y} D(\bar{x}\bar{y})$, where 
$D(\bar{x}\bar{y})$ is basic. 
Thus there is an infinite family $\Phi$ of 
basic formulas of the form $D(\bar{x}\bar{y})$
with $|\bar{x}|=k$ so that \\ 
(a) each element of $\Phi$ can be realised 
both in $(M,I^* ,\prec^* )$ and $(M,I, \prec )$, \\ 
(b) for any pair $D_1 (\bar{x}\bar{y}_1 )$ and 
$D_2 (\bar{x}\bar{y}_2 )$ from the family 
the conjunction 
$D_1 (\bar{x}\bar{y}_1 ) \wedge D_2 (\bar{x}\bar{y}_2)$
cannot be realised neither in $(M,I ,\prec )$ 
nor in $(M,I^* , \prec^* )$. 

Thus the set of $\bar{x}$-parts of realisations 
in $(M,I ,\prec )$ of formulas from $\Phi$ 
has the property contradicting to the 
property of $(I, \prec )$ from the formulation. 
%
$\Box$

\bigskip 

\begin{remark} 
{\em 
Arguments of Lemma 1.1 of \cite{iv} 
imply that if ${\bf B}_T$ from Theorem \ref{GenStr} 
satisfies Truss' condition, then there is an elementary 
embedding of $M$ into itself which takes $(I,\prec )$ into 
$(I^* ,\prec^* )$. 
} 
\end{remark}

\section{Nice enumerations}

It is known (and easily seen) that any structure having an 
AZ-enumeration is countably categorical. 
It is also clear that any AZ-enumeration is {\bf nice}, i.e. 
for any infinite $\prec$-increasing sequence $a_i \in M$ 
there are $i$ and $j$ so that the initial segment 
$\{ c: c\prec a_i \}$ is isomorphic to a substructure 
of the initial segment $\{ c: c\prec a_j \}$ by 
an automorphism mapping $a_i$ to $a_j$ \cite{az}. 

There are countably categorical structures without 
AZ-enumerations \cite{ivamaj}. 
On the other hand it is an open question if there are 
$\omega$-categorical structures without nice enumerations \cite{CaEv}. 
Since by Theorem 2.4 of \cite{CaEv} any permutation 
module of a structure having a nice enumeration, has  
the ascending chain condition for submodules it is 
also open if there are countably categorical 
structures so that there is a permutation module   
of this structure which does not have the asscending 
chain condition. 

Permutation modules which we consider appear as follows. 
For a field $F$ let $FM$ be the $F$-vector space 
where $M$ is a basis.  
Then the group algebra $FAut(M)$ naturally acts 
on $FM$, i.e. $FM$ becomes a module over 
$FAut(M)$. 
We usually consider right modules. 
When $v\in FM$, then $supp(v)$
is the set of all elements of $M$ which appear in 
$v$ with non-zero coefficients. 

\bigskip 

{\bf Some construction.} 
Let $M$ be an $\omega$-categorical 
structure. 
Consider $M$ as a relational structure 
admitting elimination of quantifiers. 
Let $I$ be an arbitrary countable subset of $M$.  
Consider a construction which can be applied  
for producing $\mathcal{L} \supset I$ with 
an AZ-well-order $<$. 
When $I$ consists of supports of elements of 
some permutation submodule the well-order 
 $<$ can be applied for finite generation of 
the submodule. 

Assume that $I$ has a well-ordering $\prec$. 
Let  $T$ be $Th(M)$ extended by all the 
sentences from $Th(M, \prec )$ of the form 
$\forall \bar{x} \neg D(\bar{x})$, where 
$D(\bar{x})$ is basic for $(M,\prec )$.
If the class ${\bf B}_T$  has the amalgamation 
property (or at least Truss' condition) one can 
apply the arguments of Theorem 1.2 of \cite{iv} 
in order to embed $(I, \prec )$ into a generic 
$(M, I^* , \prec^* )$ with countably categorical 
theory.   
In particular whenever $\bar{a}_i$, $i<\omega$, is 
a sequence of $n$-tuples from $I^*$, there exist some 
$i<j<\omega$ and an  automorphism of $(M, I^* ,\prec^* )$ 
which takes $\bar{a}_i $ to $\bar{a}_j$.  
This suggests building a required $(\mathcal{L},<)$ 
as a subset of $I^*$ (containing $I$ with $<$, 
the restriction of $\prec^*$ to $\mathcal{L}$) 
by an inductive procedure where a typical step 
looks as follows.  
We represent $I$ as a union of an infinite 
sequence of supports $\bar{b}_n$ of finite diagrams 
$D (\bar{b}_n ) \in {\bf B}_T$, $n\in \omega$. 
Assume that after step $n$ we already have 
a diagram $D (\bar{c}) \in {\bf B}_T$ where  
$\bar{b}_n \subseteq \bar{c}$.  
Amalgamating   $D (\bar{c})$  with $D (\bar{b}_{n+1})$ 
over $D (\bar{b}_n)$ we extend $<$ 
from $\bar{c}$ to some $\bar{c}'$. 
Moreover to satisfy the condition of 
AZ-orderings, in the situation when $\bar{c}_1$ and 
$\bar{c}_2$ are subtuples of $\bar{c}'$ so that 
$D (\bar{c}_1)$ is a copy of $D (\bar{c}_2)$  
under an elementary map, then amalgamating 
$D (\bar{c}')$ with itself over this map 
we obtain that $\bar{c}_1$ and $\bar{c}_2$ 
cannot be distinguished by a formula 
describing their extensions to $\bar{c}'$.  
 
There is no reason why the resulting 
$(\mathcal{L}, <)$ is a well-ordering 
(or of type $\omega$ when $I$ is so).  
Thus  we need some additional 
complications in this construction 
which would help us to control that 
$\mathcal{L}$ is decomposed into 
well-ordered parts. 
In Section 2.1 we study this issue 
in the case of well-orderings.  
In particular Lemma \ref{unions} describes 
conditions which linear orderings defined 
by $D(\bar{c})$ at our steps 
should satisfy. 
In Section 3.2 we consider a more general 
case of so called {\em quasi-well-orderings}. 
In Section 3.3 we prove a version of Theorem 2.4 
of \cite{CaEv} on permutation modules.

\subsection{Well-orderings} 

Let us recall the following notion from 
Chapter 5 of \cite{Ros}. 
\begin{quote} 
Let $(L,<)$ be a linear ordering and $x\in L$.  \\ 
Let ${\bf c} (x) = \{ y: $ the interval between $x$ and $y$ is finite $\}$. 
\end{quote} 
We consider ${\bf c}$ as a homomorphism from $(L,<)$ 
to $(L^1 ,<)$, where $L^1$ is the natural quotient of 
$L$ by the equivalence relation 
defined by the equality ${\bf c}(x) = {\bf c}(y)$.  
The ordering $L^1$ consists of intervals of $L$ and 
is called the {\bf condensation} of $L$.  
\begin{quote}  
Let ${\bf c}^0 = id$,  ${\bf c}^1 = {\bf c}$,   
${\bf c}^{\alpha +1}(x) = \{ y : {\bf c}({\bf c}^{\alpha}(x)) = {\bf c}({\bf c}^{\alpha}(y)) \}$ , \\  
and ${\bf c}^{\gamma}(x) = \bigcup \{ {\bf c}^{\alpha}(x): \alpha <\gamma \}$, where 
$\alpha$ is an ordinal and $\gamma$ is a limit ordinal.  
\end{quote} 
The appropriate condensations $L^{\alpha}$ are naturally defined. 
By Proposition 5.7 of \cite{Ros} if $x<y \in L$ and $L$ is well-ordered, 
then ${\bf c}^{\alpha} (x) = {\bf c}^{\alpha}(y)$  
if and only if $y - x <\omega^{\alpha}$. 

The least ordinal $\alpha$ such that 
${\bf c}^{\alpha} (x) = {\bf c}^{\beta}(x)$ 
for all $x\in L$ and $\beta \ge \alpha$ 
is called the ${\bf r}_F$-rank of $L$.

\begin{definicja} 
We call a finite linear ordering $L$ together 
with a function $f_c$ from the set $L\times L$ 
to the set of non-limit ordinals $On_{nl}$  
a {\bf marked ordering} if for any $\alpha \in On_{nl}$ 
the condition $f_c (x,y) \le \alpha$ defines 
an equivalence relation dividing $L$ into intervals 
and for any $\beta <\alpha$ the equivalence 
relation corresponding to $\beta$ is finer than 
the one corresponding to $\alpha$. 
We obviously assume that $f_c (x,x)=0$. 
\end{definicja} 

Note that in this situation there is an embedding 
of $L$ into a well-ordering $L'$ so that 
for any $x,y \in L$ the equality 
${\bf c}^{f_c (x,y)} (x) = {\bf c}^{f_c (x,y)} (y)$  holds in $L'$. 
Indeed let  
$0\le \gamma_1 < \gamma_2 <...<\gamma_{k-1} <\gamma_{k}$ 
be the sequence of all  values of $f_c (x,y)$. 
Note that the equivalence relation $f_c (x,y) \le \gamma_{k-1}$ 
has several, say $n>1$, equivalence classes.  
It is easy to see that $L'$ can be chosen as 
$\omega^{\gamma_k}$ where 
$\omega^{\gamma_{k-1}} \cdot n$ 
can be taken as an initial segment. 
Then each of $n$ copies of $\omega^{\gamma_{k-1}}$ 
in this segment represents  a $\gamma_{k-1}$-class. 
In each $\gamma_{k-1}$-class we find the initial 
segment consisting of appropriate copies of 
$\omega^{\gamma_{k-2}}$ and so on. 
The rest is clear. 

The following observation developes this remark. 

\begin{lem} \label{unions} 
Let $(L_i ,(f_c )_i )$ be an $\omega$-sequence 
of marked linear orderings so that 
$L_1 \subset L_2 \subset .... \subset L_i \subset ...$ 
is a sequence of ordering extensions and 
for any $i<j$ and  $x,y\in L_i$ we have 
$(f_c )_j (x,y) =  (f_c )_i (x,y)$. \\ 
For $x,y \in \bigcup L_i$  define $f_c  (x,y)$ as  $(f_c )_j (x,y)$ with  $x,y \in L_j$. \\ 
Let $\gamma_0 = lim\{ \max \{ (f_c )_i (x,y): x,y \in L_i \} :  i \rightarrow \infty \}$  
be countable. 

Assume that : \\ 
(i)  for any $a\in \bigcup L_i$ 
and any $\alpha$ of the form $f_c  (x,y)$ the set of classes 
of the $\alpha$-equivalence relation which represent in $\bigcup L_i$ 
the $(\alpha +1)$-class of $a$ is well-ordered, \\ 
(ii) there is no infinite sequence $a_1 > a_2 >... $ 
in $\bigcup L_i$ so that the sequence $f_c  (a_i ,a_{i+1} )$ 
is increasing. 

Then $\bigcup L_i$ embeds into a countable well-ordering 
$\mathcal{L}$ so that for any $x,y \in L_i \subset \mathcal{L}$ 
if the equality ${\bf c}^{\alpha} (x) = {\bf c}^{\alpha} (y)$  
holds in $\mathcal{L}$ then $f_c (x,y) \le\alpha$.   
\end{lem}

{\em Proof.}  
At every step of the construction we embed 
$L_i$ into a well-ordered set $\hat{L}_i$ 
of the form $\omega^{\gamma}$ exactly as it is described before 
the formulation of the lemma. 
Note that there is  a unique embedding of 
$\hat{L}_i$ into $\hat{L}_{i+1}$ 
so that the diagram 
constructed from $L_i \rightarrow \hat{L}_i$ 
and $L_i \rightarrow L_{i+1} \rightarrow \hat{L}_{i+1}$  
becomes commutative. 

Let $\mathcal{L}$ be the limit of  
$\{\hat{L}_i : i \rightarrow \infty\}$ with respect 
to these embeddings. 
Note that for any $\alpha \le \gamma_0$ 
each $\alpha$-equivalence class in 
$L_i$ is represented in $\hat{L}_i$ 
by a copy of  $\omega^{\alpha}$. 

To see that  $\mathcal{L}$ is well-ordered  
we prove by induction that 
\begin{quote} 
in any decreasing sequence 
$x_1 > x_2  > ... > x_k > ...$  from $\mathcal{L}$ 
for any $\alpha \le \gamma_0$, each $x_k$ 
has finitely many ${\bf c}^{\alpha}$-equivalent members. 
\end{quote} 
To see this for $\alpha \le 1$ just 
apply that each 1-class in $\bigcup L_i$ 
is well-ordered. 
Thus for any $x\in \hat{L}_i$ 
the ${\bf c}^1$-class of $x$ can not be 
extended by some elements $x'<x$ taken 
from infinitely many $\hat{L}_j$, $j>i$. 
At step $\alpha = \beta+1$
it is enough to show that any class 
of the $\alpha$-equivalence relation 
does not have an infinite decreasing 
sequence of classes of 
the $\beta$-equivalence relation 
which are represented in $\{ x_k \}$. 
Since in the extension 
$\hat{L}_i \rightarrow \hat{L}_{i+1}$ 
each new $\beta$-class in a fixed 
$\alpha$-class appears only 
together with an element of 
$L_{i+1} \setminus L_{i}$ 
this follows from the fact 
that  the set of $\beta$-classes 
of our $\alpha$-class which are 
represented in $\bigcup L_i$ is well-ordered. 

In the case of a limit $\alpha$ 
the argument is similar. 
Indeed take any $x_k$ in our 
$\alpha$-class and any 
$\beta <\alpha$.  
Now note that  in the extension 
$\hat{L}_i \rightarrow \hat{L}_{i+1}$ 
each new class for the  
$\delta$-equivalence relation with 
some $\beta\le \delta \le \alpha$
which is in  our $\alpha$-class and  less that $x_k$,  
appears only together with an element of 
$L_{i+1} \setminus L_{i}$. 
Now it suffices to apply 
conditions (i) and (ii) of the formulation. 

Since the set of ordinals is well-ordered,  
the statement we have just proved implies 
that $\mathcal{L}$ is well-ordered. 
By induction one can easily prove that for 
any $x,y \in L_i \subset \mathcal{L}$ 
if the equality 
${\bf c}^{\alpha} (x) = {\bf c}^{\alpha} (y)$  
holds in $\mathcal{L}$ then $f_c  (x,y) \le\alpha$. 
$\Box$ 

\bigskip 

The reason for the notion of marked orderings 
is as follows. 
Assume that at step $n$ of the construction 
of $(\mathcal{L}, <)$ we amalgamate   
some $D (\bar{c})$  with $D (\bar{b}_{n+1})$ 
over $D (\bar{b}_n)$ and extend $<$ 
from $\bar{c}$ to some $\bar{c}'$. 
Since $I$ is well-ordered the elements 
$\bar{b}_n$ have the structure of a marked ordering. 
If $\bar{c}$ is also a marked ordering we extend 
the corresponding function $f_c$ to 
$\bar{c}'$ so that $\bar{c}'$ becomes 
a marked ordering too. 
Our strategy is to obey the conditions 
of Lemma \ref{unions}. 
This would make $\mathcal{L}$ 
well-ordered.

\subsection{Nice quasi-well-orderings}

Let $(\mathcal{L},<)$ be a countable linearly ordered set. 
We naturally extend the betweenness realation on $\mathcal{L}$ 
to the betweenness relation on $S_1 (\mathcal{L})$, 
the set of all complete types over $\mathcal{L}$. 
In particular in the case when $p(x)$ is non-algebraic we have: 
\begin{quote} 
$c_1$ is between $p(x)$ and $c_2$ if 
$( c_2 <c_1 <x )\vee (x<c_1 <c_2 ) \in p(x)$; \\ 
$p(x)$ is between $c_1$ and $c_2$ if 
$( c_2 <x < c_1 )\vee (c_1 <x < c_2 ) \in p(x)$. 
\end{quote} 
Let $\mathcal{P}$ be a distinguished subset of 
non-algebraic types over $(\mathcal{L},<)$. 
Assume that the ordering $(\mathcal{P}, <)$ 
is discrete. 
When $p_1 <p_2$ are neighbours let 
$$ 
\mathcal{L}(p_1  ,p_2 ) = \{ x\in \mathcal{L} : x \mbox{ is between } p_1 \mbox{ and } p_2 \} \mbox{ and }  
$$ 
$$ 
\mathcal{L}^{<} (p_1 ) = \{ x \in \mathcal{L} : x < p_1  \} \mbox{ , }  
\mathcal{L}^{>} (p_2 ) = \{ x \in \mathcal{L} :  p_2  <x \} . 
$$ 
For $c_1 , c_2 \in \mathcal{L}(p_1  ,p_2 )$ 
(or  $\mathcal{L}^{<} (p_1 )$ , $\mathcal{L}^{>} (p_2 )$)                         
and $p\in \{ p_1 ,p_2 \}$  we will write $c_2 <_p c_1$ 
when $c_1$ is between $c_2$ and $p$. 
When $p_1< ...<p_n$ is a sequence from $\mathcal{P}$  
and $c_1 ,d_1 ,...,  c_n ,d_n$ belong to the corresponding 
intervals $\mathcal{L} (p_i, q_i )$ 
(or $\mathcal{L} (q_i ,p_i )$, $\mathcal{L}^{<} (p_1 )$ , $\mathcal{L}^{>} (p_n )$)   
of $\mathcal{L}$ where $q_i \in \mathcal{P}$ are appropriate neighbours, 
then we write $\bar{c} <_{\bar{p}} \bar{d}$ 
if for every $i\le n$ we have $c_i <_{p_i} d_i$.

We now generalise nice enumrations in two steps. 
At the first one we introduce 
quasi-well-orderings and then we add 
a condition which makes the ordering nice. 
Consider a structure $M$. 
Let $(\mathcal{L}, <)$ be a countable  linear order on 
a subset $\mathcal{L}\subseteq M$ and $\mathcal{P}$ 
be a finite set of  non-algebraic complete types over 
$(\mathcal{L},  <)$.  

\begin{definicja} \label{QWO} 
A linar ordering $(\mathcal{L},<)$ is 
called {\bf quasi-well-ordered} ({\bf qwo} ) 
with respect to $\mathcal{P}$ if 
$\mathcal{P}$ is presented as 
$p_1 <...< p_t$ and one of 
the following cases holds: \\ 
(a) $\mathcal{L}$ is well-ordered and $p_t$ is its $\infty$-type; \\ 
(b) $\mathcal{L}$ is well-ordered with respect to the reverse ordering and $p_1$ is its $\infty$-type; \\ 
(c) $\mathcal{L} = \mathcal{L}_1 + \mathcal{L}_2$, where $\mathcal{L}_1$ is as in (a) with the $\infty$-type, 
say $p_s$, $s\le t$,  and $\mathcal{L}_2$ is as in (b), where $p_s$ is the $\infty$-type of the reverse ordering.  
\end{definicja}

The structure 
$\omega \omega \omega^* \omega^*$ 
is a nice example of qwo with five (or three) 
non-algebraic types (here $\omega^*$ 
is a reversing ordering of $\omega$). 

Note that any subordering $\mathcal{L}(p_i , p_{i+1})$ 
is well-ordered with respect to exactly one  $<_{p_i}$ 
or $<_{p_{i+1}}$. 
This condition is very restrictive. 
The following lemma shows that a more 
general definition gives the same result.

\begin{lem} \label{GeNice} 
Assume that $(\mathcal{L},<)$ and 
$\mathcal{P} = \{ p_1 < ...  <p_t \}$ are as above. 
Then  $(\mathcal{L},<)$ is  quasi-well-ordered 
with respect to $\mathcal{P}$ if 
each of the following suborderings is 
a well-ordering : \\ 
- $\mathcal{L} (p_i ,p_{i+1} )$ with respect to $<$ or $>$, where 
$p_i  ,p_{i+1}$ is a consecutive  pair from $\mathcal{P}$, \\ 
- $(\mathcal{L}^{<} (p_1 ), <_{p_1})$ and  $(\mathcal{L}^{>} (p_t ), <_{p_t})$, 
in the case when they are not empty. 
\end{lem}  

{\em Proof.} 
The proof is easy and is based on the fact 
that each well-ordering  has the first element. 
On the other hand each element of $\mathcal{P}$ 
is non-algebraic. 
$\Box$

\bigskip

\begin{definicja}  \label{DefNice} 
Let $(\mathcal{L},<)$ be a quasi-well-ordering 
with respect to $\mathcal{P} =\{ p_1 < ...<p_k\}$. 
We say that $<$ is {\bf nice in} $M$ 
{\bf with respect to } $\mathcal{P}$ if  
for any sequence of $k$-tuples $(a^i_1 , ...a^i_k )$, 
$i<\omega$, from $\mathcal{L}$ 
(resp. $(k+1)$-tuples if $\mathcal{L}$ 
is of the form (c) of Definition \ref{QWO})  with 
$$
a^1_l \le_{p_l} ... \le_{p_l} a^i_l \le_{p_l} ... \mbox{ for } l\le k  
\mbox{ (resp. } l\le k+1 \mbox{ )  and the corresponding } 
\mathcal{L}( p_{l-1} ,p_{l} ) , 
$$ 
$$ 
\mbox{ (or }   \mathcal{L}( p_l ,p_{l+1} ) \mbox{ , } 
\mathcal{L}^{<} (p_1 ) \mbox{ , }\mathcal{L}^{>} (p_k ) \mbox{ ) } 
$$ 
there exist some $i<j<\omega$ and an automorphism  
$\xi \in Aut(M)$ such that  $\xi$ takes 
$\bar{a}_i$ to $\bar{a}_j$ and moreover for every $l\le k$ (resp. $l\le k+1$),   
$\xi$ takes the $<_{p_l}$-initial segment of $\mathcal{L}(p_{l-1} ,p_l)$ 
(resp. $\mathcal{L}(p_l ,p_{l+1})$, $\mathcal{L}^{<} (p_1 )$, $\mathcal{L}^{>} (p_k )$) 
determined by $a^i_l$  into the initial segment 
determined by  $a^j_l$.  
\end{definicja}

\begin{lem} 
Assume that $(\mathcal{L},<)$ is a countable 
quasi-well-ordering of type (a) of Definition \ref{QWO} 
with respect to $\mathcal{P} =\{ p_1 < ...<p_k\}$. 
Assume that $(\mathcal{L},<)$ is an AZ-ordering in $M$
and any elementary map between 
subsets of $\mathcal{L}$ extends to an 
automorphism of $M$. 
For example assume that $M$ is $\omega_1$-saturated. 

Then $<$ is nice with respect to $\mathcal{P}$. 
\end{lem} 

{\em Proof.} 
Take a sequence of $k$-tuples $(a^i_1 , ...a^i_k )$, 
$i<\omega$, from $\mathcal{L}$ as in 
Definition \ref{DefNice}. 
Since $<$ is an AZ-ordering there are $i<j$ 
and an elementary map which takes 
$\mathcal{L}$ into itself and takes 
$\bar{a}^i$ to $\bar{a}^j$. 
Since $a^i_l$, $l\le k$, represent each of the intervals 
$\mathcal{L}^{<} (p_1 )$ and $\mathcal{L}(p_l ,p_{l+1})$, 
$l\le k$,  
the elementary map takes the corresponding 
initial segments in appropriate way.  
$\Box$ 

\bigskip 

In the cases when $(\mathcal{L},<)$ 
is a quasi-well-ordering of type (b) or (c) 
of Definition \ref{QWO} similar 
statements hold. 
In case (b) we should only assume that 
the reverse ordering is AZ, and in case (c) 
we should assume that the ordering $\mathcal{L}_1$ 
and the reverse ordering $\mathcal{L}_2$
are AZ. 

The situation which typically arises  
in our circumstances is slightly more general. 
Let $(\mathcal{L}, <)$ be a countable linear 
order on a subset of an $\omega$-categorical 
structure $M$. 
It is usually assumed that there are 
$c_1 ,...,c_s \in \mathcal{L}$ so that each 
$\mathcal{L}(c_i ,c_{i+1})= \{ x: c_i < x < c_{i+1} \}$ 
(resp. $\mathcal{L}^< (c_1)$ and $\mathcal{L}^> (c_s )$) 
is a nice quasi-well-ordering with respect to 
some families of types $\mathcal{P}_i$ which  
represent non-algebraic types over 
$\mathcal{L}(c_i ,c_{i+1})$ 
(resp. $\mathcal{L}^< (c_1)$ and 
$\mathcal{L}^> (c_s )$).

\subsection{Permutation modules} 

Theorem 2.4 of \cite{CaEv} states that 
existence of  nice enumerations implies 
the ascending chain condition of submodules 
in permutation modules. 
It is not clear if these conditions are equivalent. 
The theorem below shows that probably a weaker 
version of nice enumerability, some kind of potential 
enumerability, still gives the same result.  
In terms of the end of the previous subsection 
it concerns the case of a single $\mathcal{L}(c_i , c_{i+1})$. 

We need the following notation. 
In the situation when  $\mathcal{L}\subseteq M$ 
and $(\mathcal{L},<)$, $\mathcal{P}$ are as above  
let $v\in FM$ and $supp(v)$ meet  
a well-ordering $\mathcal{L} (q,p)$, with $q<p\in \mathcal{P}$ 
or a reverse well-ordering $\mathcal{L} (p,q)$, 
with $p<q\in \mathcal{P}$. 
Then in the first case we denote by 
$Head^+_{p} (v)$ the maximal element 
of $supp(v)$ under the ordering $<_p$. 
In the second case we denote the corresponding element 
by  $Head^-_{p} (v)$.

\begin{thm} \label{NICE}
Let $M$ be an $\omega$-categorical 
structure and $F$ be a field. 
Assume that  $\mathcal{V} =\{ v_1 ,v_2 ,...\}$ 
is a subset of the $FAut(M)$-module 
$FM$ so that for every $i$   
$$ 
\langle \mathcal{V}\rangle_{FAut(M)} \not= \langle v_1 ,...,v_{i-1} \rangle_{FAut(M)}  . 
$$ 
Let $I$ be the union of the supports of all elements from $\mathcal{V}$.  

Then in any elementary extension of $M$ 
the set $I$ cannot be embedded into a linearly 
ordered set $(\mathcal{L},<)$  so that there is 
a family $\mathcal{P}$ such that $\mathcal{L}$ 
is a nice quasi-well-ordering  in $M$   
with respect to  $\mathcal{P}$. 
\end{thm} 

{\em Proof.} 
We start with the observation that 
the condition 
$$ 
\langle \mathcal{V}\rangle_{FAut(M)} \not= \langle v_1 ,...,v_{i-1} \rangle_{FAut(M)}  
\mbox{ for all } i \in \omega,  
$$ 
also holds in any elementary extension of $M$. 
Indeed if $M'$ is such an extension and 
$v_j \in  \langle v_1 ,...,v_{i-1} \rangle_{FAut(M)}$ 
with $i\le j$, then $v_j$ can be written as a linear combination 
of a family $u_1 ,...,u_s \in FM'$ 
where each $u_k$ is an $Aut(M')$-image of some 
of $v_1 ,...,v_{i-1}$. 
Since $M$ is $\omega$-categorical and 
$M'$ is an elementary extension of $M$,  
the type of $\bigcup_{k\le s} supp (u_k )$ over 
$\bigcup_{k<i} supp (v_k )$ is realised in $M$.   
This gives a presentation of 
$v_j$ as a linear combination of a family 
$u'_1 ,...,u'_s \in FM$ where each 
$u'_k$ is an $Aut(M)$-image of 
some  of $v_1 ,...,v_{i-1}$.   
 
Now assume that $\mathcal{L}$ contains $I$ and satisfies 
in $M$ (or in some elementary extension of $M$) 
the definition of a nice qwo with respect 
to some $\mathcal{P}= \{ p_1 < p_2 < ...< p_k\}$.    
Let $\bigcup \mathcal{L}_i$ be a partition of  
$\mathcal{L}$ into all $\mathcal{L} (p_l , p_{l+1})$ and 
$\mathcal{L}^{<} (p_{1})$,  $\mathcal{L}^{>} (p_{k})$, 
as in the corresponding definitions. 
Then $F\mathcal{L} = \oplus_i F\mathcal{L}_i$. 

Well-order $F$ in such a way that $0$ 
and $1$ are the first two elements. 
Let $<_{p_i}$ also denote the resulting lexicographic 
ordering on the corresponding  $F\mathcal{L}(p_l ,p_{l+1})$
(or for example $F\mathcal{L}^< (p_1 )$). 
It is clear that  $<_{p_i}$ is a well-ordering. 
We can now order $F\mathcal{L}$ lexicographically 
by orders $<_{p_i}$ on the corresponding components. 

We define a sequence of elements of 
$\langle \mathcal{V} \rangle_{FAut(M)} \cap F\mathcal{L}$ 
as follows. 
Let $u_1$ be the lexicographically minimal element of 
$(\langle \mathcal{V}\rangle_{FAut(M)} \cap F\mathcal{L}) \setminus \{ 0\}$. 
For each $s$ let $u_s$ be the lexigraphically 
minimal element of 
$$ 
(\langle \mathcal{V}\rangle_{FAut(M)} \cap F\mathcal{L}) \setminus  \langle u_1 ,...,u_{s-1} \rangle_{FAut(M)} .
$$ 

Fix $<_{p_l}$ and $F\mathcal{L}(p_{l-1} ,p_{l})$ 
(or $F\mathcal{L}(p_l,p_{l+1} )$, but we assume the  $+$-case) 
with the minimal number in the decomposition  
$F\mathcal{L} = \oplus_i F\mathcal{L}_i$  
so that there is an infinite subsequence of $u_1 ,u_2 ,...$ 
with $<_{p_l}$-increasing 
$F\mathcal{L}(p_{l-1} ,p_l )$-components, say 
$$
u_{1l} <_{p_l} u_{2l}  <_{p_l} u_{3l} <_{p_l} ... <_{p_l} u_{sl} <_{p_l} ... 
\mbox{ . }
$$ 
We may assume that $u_1 , u_2 , ...$ have 
the same components enumerated before 
$F\mathcal{L}(p_{l-1} ,p_l)$ 
(then such components are equal to $0$) 
and $Head^+ (u_{sl} )$   
occurs in $u_s$ with coeffecient $1$. 
Moreover we may assume that for all $s$, 
$Head^+ (u_{sl} )< Head^+ (u_{(s+1)l})$. 
Indeed, if $Head^+ (u_{sl} )= Head^+ (u_{(s+1)l})$, 
then the element $u_{s+1} - u_{s}$ 
contradicts the choice of $u_{s+1}$ as 
the minimal element outside 
$\langle u_1 ,...,u_{s-1} \rangle_{FAut(M)}$. 

So $\{ Head^+ (u_{sl} ) :si\in \omega \}$
is an infinite subset of $\mathcal{L}(p_{l-1} ,p_l )$. 
Since  $\mathcal{L}$ is a nice qwo in $M$ 
with respect to $\{ p_1 ,...,p_k \}$
there are numbers $s<t$  
and an automorphism  $\xi \in Aut(M)$ 
which takes $Head^+ (u_{sl} )$ to  $Head^+ (u_{tl} )$  
and the initial segment of $<_{p_l}$ determined by $u_{sl}$
to the initial segment determined by $u_{tl}$.  
Moreover we may assume that 
$supp ((u_s )\xi ) \subset \mathcal{L}$. 
Then $u_t - (u_s )\xi $ is less that $u_t$ 
under the lexicographic ordering. 
This  contradicts the choice of $u_t$. 
$\Box$ 

\bigskip 

There are some interestion special cases. 
For example assume that $(I, < )$ 
is qwo, but additionally $(I,< )$ is 
order indiscernible in $M$.  
By $\omega$-categoricity of $M$ we see that:  
\begin{quote} 
if any order-preserving injective map $I \rightarrow I$ 
extends to an automorphism of $M$, then 
$(I,<)$ is a nice qwo in $M$. 
\end{quote} 
This in particular holds when $M$ is $\omega_1$-saturated. 
By Theorem \ref{NICE} the set $I$ cannot be the union of 
supports of a subset $\mathcal{V}$ of a permutation module 
$FM$ generating a submodule which is not finitely generated.

\subsection{Summary and final remarks} 

We now summarise the approach of 
the paper. 
Let $M$ be an $\omega$-categorical 
structure. 
Consider $M$ as a relational structure 
admitting elimination of quantifiers. 
It is an open question if $M$ has a nice enumeration. 
We try to replace this property by 
the following approach.  
Let $I$ be an arbitrary countable subset of $M$.  
Apply the construction which we described in 
the beginning of Section 3.1. 
Our goal is some $(\mathcal{L},\prec)$ contradicting 
the conclusion of Theorem \ref{NICE} 
(when we have such $\mathcal{L}$ we know that 
$I$ is not the union of supports of a sequence 
generating a permutation submodule which is not 
finitely generated). 
As above we assume that $I$ is well-ordered 
by some $\prec$. 
Let $T$ be the theory as in the formulation of 
Theorem \ref{GenStr}.  
By Lemma \ref{HuNe} if the class ${\bf B}_T$ has 
the Ramsey property, then it has the amalgamation 
property. 
Then we represent $I$ as a sum of an infinite sequence 
of finite diagrams $D (\bar{b}_n ) \in {\bf B}_T$, $n\in \omega$. 
If after step $n$ we already have 
a diagram $D (L_n ) \in {\bf B}_T$ where  
$\bar{b}_n \subseteq L_n$, then 
amalgamating   $D (L_n )$  with $D (\bar{b}_{n+1})$ 
over $D (\bar{b}_n)$ we extend $\prec$ 
from $L_n$ to some $L_{n+1}$. 

We build $\mathcal{L} = \bigcup L_n$ 
so that $\mathcal{L}$ can be presented as 
a finite union of substructures of the 
form $\mathcal{L}(c_i ,c_{i+1})$  
as in the end of Section 3.2. 
Our strategy is to arrange that each 
sequence 
$\mathcal{L}(c_i ,c_{i+1})\cap L_n$, 
$n \in \omega$,   
satisfies the conditions of Lemma \ref{unions}  

\bigskip

There are some additional 
issues when 
$\mathcal{L}$ can be also ordered, say by $<$, 
so that $(\mathcal{L} ,<)$ is order 
indiscernible in $M$. 
The case when $\mathcal{L}$ is an idiscernible 
set with respect to $Th(M)$ in fact 
appeared in the remark after 
Theorem \ref{NICE} with $\mathcal{L} = I$. 
It is worth noting that in this case 
the constructed 
$\prec$ is an indiscernible ordering in $M$. 
Moreover applying $\omega$-categoricity 
of $M$ we see that it is an AZ-ordering 
in sufficiently saturated extensions of $M$.  

It is also worth mentioning that 
in this case if additionally $(\mathcal{L}, \prec )$ 
is a dense ordering, then taking any $<$ 
of type $\omega$ we obtain a nice enumeration 
of the structure $(\mathcal{L}, \prec )$. 
To see this one can apply an observation 
from \cite{CaEv} that any $\omega$-enumeration 
of a dense linear ordering is nice. 
 
Consider the opposite case, i.e. assume that 
$\mathcal{L}$ is an order indiscernible sequence 
of type $\omega$ with respect to $Th(M)$,   where 
$<$ is the corresponding ordering of $\mathcal{L}$,  
and the type (with respect to $Th(M)$) of any 
$<$-increasing tuple from $\mathcal{L}$ differs from 
the type of any non-trivial permutation of the tuple.  
Let $T$ be as above and let $\mathcal{C}$ be 
the subclass of all members of ${\bf B}_T$ of the 
form $D(\bar{a})$ with $\bar{a}\subset \mathcal{L}$. 
By Section 3 of \cite{AlCh} 
\begin{quote} 
if all possible finite linear orderings can be 
realised by $\prec$ on $<$-increasing 
tuples from $\mathcal{L}$,   then  
the family $\mathcal{C}$ has an infinite 
anti-chain with respect to the natural 
embedding defined in ${\bf B}_T$.    
\end{quote}
By \cite{AlCh} this condition means that there are 
enumerations of $\mathcal{L}$ which are not AZ-enumerations. 

\begin{remark} \label{RemRam} 
{\em 
Assume that $(\mathcal{L}, \prec )$ is as above  
and  $I$ is order indiscernible 
in $M$ with respect to some order $<$. 
Then the subclass of all members of $\mathcal{C}$  
has the Ramsey property with respect to $<$-preserving 
embeddings if and only if so does 
the class of all $<$-ordered finite structures 
of the relation $\prec$ realisable in $(\mathcal{L},<)$. 
}
\end{remark} 

\begin{remark} 
{\em 
Let us notice that if a countable structure $M$ is a model  
of an $\omega$-categorical universal theory, 
then any enumeration of $M$ is nice.  
This follows from Theorem 1 of \cite{palyutin} stating 
the existence of a function $s:\omega \rightarrow \omega$, 
so that for any substructure $B<M$ of size $\ge s(n)$ 
any 1-type over an $n$-element subset of $B$ 
is realised in $B$ (i.e. in particular any 
$n$-type over any $b\in B$ is realised in $B$). 
Now let $p(x)$ be a non-algebraic type over $\emptyset$. 
If $a_1 ,...,a_n ,...$ is an infinite sequence of realisations 
of $p(x)$ in order of the enumeration of $M$, then 
choosing $n$ sufficiently big, we find a realisation 
of the initial segment defined by $a_1$ in the initial 
segment of $a_n$, where $a_n$ plays the role of $a_1$. 
Now the condition of nice enumeration can be easily verified.    

We in particular have that permutation modules over 
$M$ satisfy the ascending chain condition. 
Let us also mention an old open question if there is 
an $\omega$-categorical universal theory which is not 
$\omega_1$-categorical \cite{palyutin}. 
}
\end{remark}

\bigskip 

{\bf Anti-chains of substructures.} 
We now prove a proposition which  
connects the ascending chain condition 
for submodules of permutation modules 
with the existence of infinite anti-chains 
of substructures. 

\begin{prop} \label{RED} 
Let $M$ be an $\omega$-categorical 
structure. 
Considering $M$ as a relational structure assume 
that the the family $\mathcal{K}$ of all finite substructures  
of $M$ does not have an infinite antichain with 
respect to embedding induced by automorphisms of $M$. 
 
Then for any finite field $GF(q)$ the  permutation module 
$GF(q) M$ over $GF(q) Aut(M)$ has the ascending 
chain condition for submodules. 
\end{prop} 

It is worth also noting that the proposition in fact 
reduces the case of the ascending chain condition  
of permutation modules over finite fields to the 
case of permutation modules over $GF(2)$  
(the permutation module 
$GF(2)M$ over $GF(2) Aut(M)$ can be identified 
with $\mathcal{K}$). 

\bigskip

{\em Proof of Proposition \ref{RED}.} 
Since $M$ is countable, it can be eumerated: 
$M= \{ a_1 ,a_2 ,...\}$. 
Suppose that there is a sequence $p_i \in GF(q)M$, 
$i\in\omega$,  which generates a non-finitely generated 
$GF(q)$-permutation submodule of $GF(q)M$. 
We express each $p_i$ as a sum 
$$ 
f_1 (\sum_{j\in D_1} a_j )+  f_2 (\sum_{j\in D_2} a_j )+ ...+  f_{q-1} (\sum_{j\in D_{q-1}} a_j ) , 
$$ 
where $GF(q)\setminus \{ 0\}$ is enumerated as 
$\{ f_1 ,...,f_{q-1}\}$ and sets $D_k$ are pairwise disjoint.   
We may assume that for sufficiently large $i$ 
the element $p_i$ always has the minimal number of 
non-zero sums $\sum_{j\in D_k} a_j$ among all 
elements from the difference 
$$
 \langle p_0 ,..,p_{i} \rangle_{GF(q)} \setminus   \langle p_0 ,..,p_{i-1}\rangle_{GF(q)} .
$$ 
Then there is an infinite subsequence of $p_i$ 
consisting of elements which have non-zero 
members $f_k (\sum_{j\in D_k} a_j )$ for the same $k$. 
Moreover we may assume that any $p_i$ of this 
subsequence satisfies the condition 
$$ 
1\le |D_1 |\le min(|D_1 |,|D_2 |,...,|D_{l-1} | ) \mbox{ with } 
D_l =...=D_{q-1}=\emptyset . 
$$ 
We may also assume that the size $|D_1 |$ in these $p_i$ is minimal
among all elements from the difference 
$$
\langle p_0 ,..,p_{i} \rangle_{GF(q)} \setminus   \langle p_0 ,..,p_{i-1}\rangle_{GF(q)} 
$$ 
with $D_l =...=D_{q-1}=\emptyset$. 

We now identify each $D_1$ with a substructure from $\mathcal{K}$. 
Using the assumption of the theorem  we choose two members 
$p_i =f_1 (\sum_{j\in D^i_1}a_j ) + ...$ and $p_l=f_1 (\sum_{j\in D^l_1}a_j ) + ...$ 
of the subsequence described above so that 
$D^i_1$ embeds into $D^l_1$ by a map extending to 
an automorphism of $M$, say $\alpha$. 
Then the element $p_l -\alpha p_i$ contradicts the choice of $p_l$.  
$\Box$ 



\bigskip

Institute of Mathematics, University of Wroc{\l}aw, pl.Grunwaldzki 2/4, 50-384 Wroc{\l}aw, Poland, \\
 E-mail: ivanov@math.uni.wroc.pl 

\end{document}